\documentstyle[12pt]{article}

\newtheorem{thm}{Theorem}
\newtheorem{lem}[thm]{Lemma}

\newtheorem{prop}[thm]{Proposition}
\newenvironment{pf}{\noindent {\em Proof:}}{$\Box$\\}

\newcommand{\N}{\mbox{\hskip.1em N \hskip -1.25em \relax I
\hskip .1em}}
\newcommand{\R}{\mbox{\hskip.1em R \hskip -1.25em \relax I
\hskip .1em}}
\newcommand{\seq}{\el^{p,\infty}}
\newcommand{\unit}{L^{p,\infty}[0,1]}
\newcommand{\rline}{L^{p,\infty}[0,\infty)}
\newcommand{\ok}{L^{p,\infty}[0,k]}
\newcommand{\rnorm}{\biggr\|}
\newcommand{\lnorm}{\biggl\|}
\newcommand{\llang}{\biggl\langle}
\newcommand{\lrang}{\biggr\rangle}
\newcommand{\lv}{\biggl|}
\newcommand{\rv}{\biggr|}
\newcommand{\lp}{\bigg(}
\newcommand{\rp}{\bigg)}
\newcommand{\msp}{(\Omega ,\Sigma ,\mu)}
\newcommand{\cembeds}{\stackrel{c}{\hookrightarrow}}
\newcommand{\el}{\ell}
\newcommand{\ktwon}{k\!\cdot\! 2^n}
\newcommand{\wktwon}{k\cdot 2^n}
\newcommand{\ktwom}{k\!\cdot\! 2^m}
\newcommand{\intf}{\int^{j/2^n}_{(j-1)/2^n} f}
\newcommand{\ra}{\Rightarrow}
\newcommand{\phinj}{\Phi_{n,j}}
\newcommand{\phin}{\Phi_{n,1}}
\newcommand{\sumj}{\sum^{\wktwon}_{j=1}}
\newcommand{\sumn}{\sum^{l-1}_{n=0}}
\newcommand{\sumjn}{\sum^{2^n}_{j=1}}

\begin{document}

\begin{center}
  {\Large Isomorphism of certain weak $L^p$ spaces}\vspace{3mm}\\
  {\large Denny H. Leung}
\end{center}

\vspace{1mm}

\begin{abstract}
   It is shown that the weak $L^p$ spaces $\el^{p,\infty},
L^{p,\infty}[0,1]$, and $L^{p,\infty}[0,\infty)$ are isomorphic as
Banach spaces.
\end{abstract}

\begin{figure}[b]
  \rule{3in}{.005in}
\end{figure}

\baselineskip 4ex

\section{Introduction}

The Lorentz spaces play an important role in interpolation theory.
They also form a class of Banach spaces generalizing the classical
$L^p$ spaces.  In this paper, we continue the comparison of the Banach
space structures of various weak $L^p$ spaces begun in \cite{L1} and
\cite{L2}.  In \cite{L1}, mimicking the construction of the Rademacher
functions, it was shown that $\el^2$ can be embedded complementably
into $\seq$.  In \cite{L2}, we showed that $\seq$ can in turn be
embedded complementably (even as a sublattice) into $\unit$.  Here, we
complete and extend these results by showing that, in fact, the three
weak $L^p$ spaces $\seq , \unit$, and $\rline$ are isomorphic as Banach
spaces. This question was also mentioned in \cite{CD}.\\
\indent We start by recalling some standard definitions.  Let $\msp$
be an
arbitrary measure space.  For $1 < p < \infty$, the weak $L^p$ space
$L^{p,\infty}\msp$ is the space of all $\Sigma$-measurable functions
$f$ such that $\{\omega : |f(\omega)| > 0\}$ is $\sigma$-finite and
\[ \|f\| \equiv \sup_B \frac{\int_B |f| d\mu}{\mu (B)^{1/q}} < \infty
,\]
where $q = p/(p-1)$, and the supremum is taken over all measurable
sets $B$ with $0 < \mu(B) < \infty$.  If $\msp$ is a real interval $I$
endowed with Lebesgue measure, we write $L^{p,\infty}(I)$; while
$\el^{p,\infty}$ and $\el^{p,\infty}(m)$ will stand for the weak $L^p$
spaces on $\N$ and $\{1,2,\ldots,m\}$ respectively, both with the
counting measure.\\
\indent For a real valued function defined on $\msp$, let $f^*$ denote
the
decreasing rearrangement of $|f|$ \cite{LT2}; similarly for $(a^*_n)$,
where $(a_n)$ is a real sequence.  It is well known that the
expression
\[ |||f||| \equiv \sup_{t>0}\, t^{1/p}f^*(t) \]
satisfies $|||f||| \leq \|f\| \leq q|||f|||$.  Let $L^{q,1}\msp$
denote the space of all measurable functions $f$
such that $\|f\|_{q,1} = \int^\infty_0 t^{-1/p}f^*(t) dt < \infty$.
Then $L^{p,\infty}\msp$ is naturally isomorphic to the dual of
$L^{q,1}\msp$, where the isomorphism constant depends only on $p$.
Finally, we note that any weak $L^p$ space satisfies
an upper $p$-estimate with constant $1$ \cite{LT2}.

\section{Main Theorem} \label{theorem}

Fix $1<p<\infty$ and $q = p/(p-1)$, our goal is to prove the following

\begin{thm} \label{main}
The Banach spaces $\seq, \unit $, and $\rline $ are isomorphic.
\end{thm}

Since the proof of the Theorem goes through several intermediate
embeddings and is rather circuitous, we first give an outline of the
procedure.  We start by recalling the following well known variant of
Pe\l czynski's ``decomposition method'' \cite{LT1}.  The {\em square}
of a Banach space $E$ is the Banach space $E\oplus E$.

\begin{thm}
Let $E, F$ be Banach spaces which are isomorphic to their squares.
Suppose that each is isomorphic to a complemented subspace of the
other, then they are isomorphic.
\end{thm}

\begin{pf}
Using the symbol ``$\sim$'' for ``is isomorphic to'', we find Banach
spaces $G$ and
$H$ such that $E \sim F \oplus G$ and $F \sim E \oplus H$.  Then
\[ E \oplus F \sim E \oplus (E \oplus H) \sim E \oplus H \sim F
,\hspace{2em} {\rm and} \]
\[ E \oplus F \sim (G \oplus F) \oplus F \sim G \oplus F \sim E .\]
\end{pf}

It has been mentioned already that $\seq$ embeds complementably into
$\unit$ \cite{L2}, and it is clear that $\unit$ embeds complementably
into $\rline$.  Since these spaces are obviously isomorphic to their
squares, the proof of Theorem \ref{main} will be complete if we can
show that $\rline$ embeds complementably into $\seq$.  This we do in a
number of steps.\\
\indent For any $k \in \N$, let $X_k = \biggl(\sum^\infty_{n=0}\oplus
\seq
(\ktwon)\biggr)_{\el^\infty}$. Then, using ``$\cembeds$'' to denote
``embeds complementably into'', we will show that
\begin{eqnarray*}
 \rline & \cembeds & \biggl(\sum^\infty_{k=1}\oplus
L^{p,\infty}[0,k]\biggr)_{\el^\infty} \\
  & \cembeds & \biggl(\sum^\infty_{k=1}\oplus X_k\biggr)_{\el^\infty}\\
  & \cembeds & X_1 \\
  & \cembeds & \seq .
\end{eqnarray*}
The third embedding in this chain is obvious, and the first is also
straight forward.  The second embedding is accomplished by showing
that $L^{p,\infty}[0,k] \cembeds X_k$ with uniform constants.  This
relies on the techniques of \cite{L1}.  For the last link in the
chain, we show that $X_1$ is isomorphic to a weak$^*$ closed subspace
of
$\seq$ generated by long blocks with constant coefficients.  The
complementation is then effected by a conditional expectation operator
\cite{LT2}.

\section{Finding $\ok$ in $X_k$}

\begin{lem} \label{boundf}
Let $k \in{\rm\N\ }$and let $f \in \rline$.  Then
\[ \sup_n 2^{n/q} \lv\lv\lv\lp\int^{j/2^n}_{(j-1)/2^n}
f\rp^{\wktwon}_{j=1}\rv\rv\rv \leq \|f\| .\]
\end{lem}

\begin{pf}
Fix $n$.  Since the expressions involved are rearrangement invariant,
we may assume that $|\int^{j/2^n}_{(j-1)/2^n} f|$ decreases with $j$.
Let $a_j = j^{1/p} \lv\int^{j/2^n}_{(j-1)/2^n} f\rv$.  Then for
 $1 \leq
j \leq k\cdot 2^n$,
\begin{eqnarray*}
\int^{j/2^n}_{0} |f| & \geq & \sum^{j}_{i=1}\lv\int^{i/2^n}_{(i-1)/2^n}
f\rv \\
& \geq & j^{1/q} a_j.
\end{eqnarray*}
Hence
\begin{eqnarray*}
\|f\| & \geq & (j/2^n)^{-1/q}\int^{j/2^n}_0 |f| \\
& \geq & 2^{n/q} a_j .
\end{eqnarray*}
Taking the supremum over $j$ finishes the proof.
\end{pf}

We introduce some more notation.  An element $x \in X_k$ will be
written as $x = (x_n)^\infty_{n=0}$, where each $x_n \in \seq
(\wktwon)$. Each $x_n$  is in turn a finite real sequence
$(x_n(j))^{k\cdot 2^n}_{j=1}$. Define
\[ Y_k = \{x \in X_k : x_n(j) =
2^{-1/q}(x_{n+1}(2j-1)+x_{n+1}(2j)),\hspace{1em} 1
\leq j \leq k\cdot 2^n,\hspace{1em} n \geq 0 \}.  \]
We will say that a linear operator $T$ mapping between Banach spaces
is a
$K$-iso\-mor\-phi\-sm (into) if $K^{-1}\|x\| \leq \|Tx\| \leq K\|x\|$.

\begin{prop} \label{embed}
Define $T_k : L^{p,\infty}[0,k] \to X_k$ by $T_kf = x$, with
$x_n(j) = 2^{n/q} \intf$ for $1 \leq j \leq \ktwon, n \geq 0$.  Then
$T_k$ is
a $q-$isomorphism of $\ok$ onto $Y_k$.
\end{prop}

\begin{pf}
Fix $f \in \ok$ and let $x = T_kf$.  By Lemma \ref{boundf},
$\sup_n|||x_n|||
\leq \|f\|$.   Thus
\[ \|T_kf\| = \sup_n\|x_n\| \leq q\sup_n|||x_n||| \leq q\|f\|. \]
On the other hand, letting $E_n$ denote the conditional expectation
operator \cite{LT2} with respect to the partition $([(j-1)/2^n,
j/2^n])^{\wktwon}_{j=1}$, it is easy to see that $|||x_n||| =
|||E_nf|||$. Now since $E_nf \to f$ in the weak$^*$ topology, we have
\begin{eqnarray*}
\|f\| & \leq & \limsup \|E_nf\| \\
 &  \leq  & q\limsup |||E_nf||| \\
 & = & q\limsup |||x_n||| \\
 & \leq & q\limsup \|x_n\| \\
 & \leq & q\|T_kf\|
\end{eqnarray*}
This proves that $T_k$ is a $q-$isomorphism.  Clearly, $T_k$ maps into
 $Y_k$.
Conversely, for $x \in Y_k$, let $f_n = 2^{n/p} \sum^{\wktwon}_{j=1}
x_n(j)\chi_{n,j}$, where $\chi_{n,j}$ is the characteristic function of
the interval $[(j-1)/2^n,j/2^n]$.  Then $|||f_n||| = |||x_n|||$.  In
particular, $(f_n)$ is a bounded sequence in $\ok$.  Now if $f$ is a
weak$^*$ cluster point of the sequence $(f_n)$, then, using the fact
that $x \in Y_k$,  it is easy to see
that $T_kf = x$.  Hence $T_k$ maps onto $Y_k$, as required.
\end{pf}

We proceed to show that $Y_k$ is complemented in $X_k$.  Fix $k \in
\N$, let
\[ Z_k = \{x \in X_k : \mbox{There exists}\ i\ \mbox{such that\ }
\sum^{\wktwon}_{j=1}x_n(j) = 0\ \mbox{for all\ } n \geq i\}. \]
It is clear that $Z_k$ is a linear subspace of $X_k$.  Also define $u
\in
X_k$ so that $u_n(j) = (\ktwon)^{-1/p}, \hspace{.5em}1 \leq j \leq
\ktwon, \hspace{.5em} n \geq 0$.

\begin{lem} \label{boundedness}
Define $\phi : {\rm span}\{Z_k,\{u\}\} \to {\rm\R\ }$by
\[ \phi(z+au) = a\  \mbox{for all\ } z \in Z_k\ \mbox{and\ } a \in\
{\rm\R\,}.  \]
Then $\|\phi\| \leq 1$ with respect to the norm on $X_k$.
\end{lem}

\begin{pf}
For all $n \geq 0$, the functional $x'_n$ on $X_k$ given by $x'_n(x) =
\sum^{\wktwon}_{j=1}x_n(j)$ has norm $\leq (\ktwon)^{1/q}$.  Now $z
\in
Z_k$ implies $x'_n(z) = 0$ for all large $n$.  Thus
\[
\begin{array}{lcr}
\hspace{10em} & \|z+au\| \leq 1 & \hspace{10em} \vspace{1ex}\\
\ra & |x'_n(z+au)| \leq (\ktwon)^{1/q} & \vspace{1ex}\\
\ra & |a|x'_n(u) \leq (\ktwon)^{1/q}  & \vspace{1ex}\\
\ra & |a| \leq 1. &
\end{array} \]
\end{pf}

Let $\Phi$ be a norm
preserving extension of $\phi$ to all of $X_k$.
For an element $x \in X_k$, the {\em support} of $x$, supp\,$x =
\{(n,j): x_n(j) \neq 0\}$. As it is clear that any finitely supported
element of $X_k$ is in $Z_k$, we must have $\Phi = 0$ on
$\lp\sum^\infty_{n=0}\oplus \seq (\ktwon)\biggr)_{c_0}$.\\
\indent If $A$ is a subset of $\Gamma = \{(n,j): 1 \leq j
\leq \ktwon, n \geq 0 \}$, The operator on $X_k$ given by
multiplication by the characteristic function of $A$, which we denote
by $\chi_A$, is a norm one projection on $X_k$. For any $(n,j) \in
\Gamma$, let $A_{n,j} =$ supp $T_k\chi_{[(j-1)/2^n,j/2^n]}$.  Then let
$\Phi_{n,j} = \chi'_{A_{n,j}}\Phi$.

\begin{lem} \label{sum}
For every $n \geq 0$, $(\Phi_{n,j})^{\wktwon}_{j=1}$ is a sequence of
pairwise disjoint (in the lattice sense) functionals on $X_k$.
Moreover, $\Phi_{n+1,2j-1}+\Phi_{n+1,2j} = \Phi_{n,j}$ for all $ n \geq
0,\ 1\leq j\leq \ktwon$.
\end{lem}

\begin{pf}
Note that $A_{n+1,2j-1}\cup A_{n+1,2j} = A_{n,j}$.  Then both
assertions follow easily from the fact that $A_{n,j}\cap
A_{n,j'}$ is finite for $1\leq j\neq j' \leq \ktwon$, and $\Phi$ is
$0$ on finitely supported elements.
\end{pf}

Given $n\geq 0,\ 1\leq j\leq \ktwon$, define $S_{n,j}:X_k \to X_k$ by
$S_{n,j}x = z$, where $z_m = 0$ for $0\leq m < n$, and $z_m(i) =
x_m(i+(j-1)2^{m-n}\ \mbox{mod}(\ktwom))$ for $m\geq n,\ 1\leq i \leq
\ktwom$.  It is clear that $S_{n,j}$ has norm $1$, and that
$x-S_{n,j}x \in Z_k$.\\

\begin{lem} \label{norm}
For $n\geq 0,\ 1\leq j\leq \ktwon,\ \|\Phi_{n,j}\| = \|\Phi_{n,1}\|$.
\end{lem}

\begin{pf}
For any $x \in X_k$,
\[ \phinj(x) = \Phi(\chi_A{_{n,j}}x) = \Phi(S_{n,j}\chi_{A_{n,j}}x) \]
since $\Phi = 0$ on $Z_k$.  But the support of
$S_{n,j}\chi_{A_{n,j}}x$ is contained in $A_{n,1}$.  Therefore,
\[ \phinj(x) = \Phi_{n,1}(S_{n,j}\chi_{A_{n,j}}x) \leq \|\Phi_{n,1}\|
\|x\| \]
since both $S_{n,j}$ and $\chi_{A_{n,j}}$ are norm $1$ operators.
Thus $\|\phinj\| \leq \|\Phi_{n,1}\|$.  The reverse inequality can be
obtained similarly.
\end{pf}

\begin{lem} \label{normphi}
For all $n \geq 0$, $\|\phin\| \leq (\ktwon)^{-1/q}\|\Phi\|$.
\end{lem}

\begin{pf}
Since $\seq$ satisfies an upper $p$-estimate with constant $1$, so
does $X_k$.  Hence $X'_k$ satisfies a lower $q$-estimate with
constant $1$. By Lemma \ref{sum}, $\Phi = \sum^{\wktwon}_{j=1}\phinj$
and the summands are pairwise disjoint.  Then by Lemma \ref{norm},
\[ \|\Phi\| \geq \lp\sum^{\wktwon}_{j=1}\|\phinj\|^q\rp^{1/q} =
(\ktwon)^{1/q}\|\phin\| . \]
\end{pf}

\begin{lem} \label{bound}
Given real numbers $b_1, b_2, \ldots, b_i$ and a sequence $(c_j)$ in
the unit ball of $\seq$,
\[ \lnorm\lp\sum^i_{j=1}b_jc_{li+j}\rp^{\infty}_{l=0}\rnorm_{\seq}
    \leq
     q^2\sum^i_{j=1}b^*_j(j^{1/q}-(j-1)^{1/q}) .\]
\end{lem}

\begin{pf}
Let $K = \|(\sum^i_{j=1}b_jc_{li+j})^{\infty}_{l=0}\|_{\seq}$.
By rearranging $(c_j)$, we may assume without loss of generality that
$|\sum^i_{j=1}b_jc_{li+j}|$ is a decreasing function of $l$.  Recall
that $|||x||| \geq \|x\|/q$ for all $x \in \seq$. Thus, given
 $\epsilon >
0$, there exists
$r$ such that $r^{1/p}|\sum^i_{j=1}b_jc_{ri+j}| > (K-\epsilon)/q$.
Then, since $c^*_j \leq j^{-1/p}$ for all $j \geq 1$,
\begin{eqnarray*}
\frac{(K-\epsilon)(r+1)^{1/q}}{q} & \leq &
        \sum^r_{l=0}\lv\sum^i_{j=1}b_jc_{li+j}\rv \\
 & \leq & \sum^i_{j=1}\sum^r_{l=0}|b_jc_{li+j}| \\
 & \leq & \sum^i_{j=1}b^*_j\sum^{r}_{l=0}((j-1)(r+1)+l+1)^{-1/p} \\
 & \leq & q(r+1)^{1/q}\sum^i_{j=1}b^*_j(j^{1/q}-(j-1)^{1/q}) .
\end{eqnarray*}
Multiplying by $q(r+1)^{-1/q}$ finishes the proof.
\end{pf}

\begin{lem} \label{compare}
For any real sequence $(a_j)$,
\[ \lnorm\sum^{\wktwon}_{j=1}a_j\phinj\rnorm \leq qk^{-1/q}
    \lnorm\sum^{\wktwon}_{j=1}a_j\chi_{[(j-1)/2^n,j/2^n]}\rnorm_{q,1}\
, \]
where $\|\cdot\|_{q,1}$ denotes the norm on $L^{q,1}[0,k]$.
\end{lem}

\begin{pf}
Fix $x \in X_k$ with norm $\leq 1$,
\begin{eqnarray*}
\llang x, \sum^{\wktwon}_{j=1}a_j\phinj\lrang & = &
      \sumj a_j\langle\chi_{A_{n,j}}x, \Phi\rangle \\
 & = & \sumj a_j\langle S_{n,j}\chi_{A_{n,j}}x, \Phi\rangle \\
 & = & \langle y, \phin\rangle ,
\end{eqnarray*}
where $y = \sumj a_jS_{n,j}\chi_{A_{n,j}}x$.  By Lemma \ref{bound},
$\|y_m\| \leq q^2\sum^{\wktwon}_{j=1}a^*_j(j^{1/q}-(j-1)^{1/q})$ for
 all $m
\geq 0$.  Hence
\begin{eqnarray*}
\llang x, \sum^{\wktwon}_{j=1}a_j\phinj\lrang & \leq & \|y\| \|\phin\|
\\
 & \leq &
      q^2\|\Phi\|\sum^{\wktwon}_{j=1}a^*_j
             \frac{j^{1/q}-(j-1)^{1/q}}{(\ktwon)^{1/q}} \\
 & \leq & qk^{-1/q}
      \lnorm\sum^{\wktwon}_{j=1}a_j\chi_{[(j-1)/2^n,j/2^n]}\rnorm_{q,1}
   .
\end{eqnarray*}
Here the second inequality follows form Lemma \ref{normphi}, while the
 third
inequality is true because $\|\Phi\| \leq 1$.
\end{pf}

\begin{prop} \label{eachk}
Define $P_k:X_k\to X_k$ by $P_kx=y$, where
\[ y_n(j) = (\ktwon)^{1/q}\langle x,\phinj\rangle \]
for all $n\geq 0,\ 1 \leq j \leq \ktwon$. Then $P_k$ is a projection
from $X_k$ onto $Y_k$ of norm $\leq q^2$.
\end{prop}

\begin{pf}
Since $(\ell^{q,1})' = \seq$ isomorphically, there is a constant $K$
(actually $q$ suffices) such that
\[ \|(a_j)\| = K\sup_{(b_j) \in U}\sum a_jb_j \]
for all $(a_j) \in \seq$, where $U$ denotes the unit ball of
$\ell^{q,1}$.  Thus, given  $x \in X_k$ and $n \geq 0$,
\begin{eqnarray*}
(\ktwon)^{1/q}\|(\langle x,\phinj\rangle)^{\wktwon}_{j=1}\| & \leq &
  q(\ktwon)^{1/q}\sup_{(b_j)\in U}\sum b_j\langle x,\phinj\rangle \\
 & \leq & q(\ktwon)^{1/q}\sup_{(b_j)\in U}\|x\|\|\!\sum b_j\phinj\| \\
 & \leq & q^22^{n/q}\|x\|\sup_{(b_j)\in U}
   \lnorm\sum^{\wktwon}_{j=1}b_j\chi_{[(j-1)/2^n,j/2^n]}\rnorm_{q,1}
     \hspace{1em} \mbox{by Lemma \ref{compare}} \\
 & = & q^2\|x\| .
\end{eqnarray*}
Hence $\|P_k\| \leq q^2$.  Using Lemma \ref{sum}, it is easy to see
that $P_k$ maps into $Y_k$. Conversely, let $y \in Y_k$.
Then for $ n\geq 0$, $1 \leq j \leq \ktwon$,
\[ y_n(j)(\ktwon)^{1/p}\chi_{A_{n,j}}u - \chi_{A_{n,j}}y \in Z_k. \]
 Note
also that $S_{n,j}\chi_{A_{n,j}}u - \chi_{A_{n,1}}u \in Z_k$.  Hence
for $1 \leq j \leq \ktwon$,
\[ \Phi(\chi_{A_{n,j}}u) = \Phi(S_{n,j}\chi_{A_{n,j}}u) =
   \Phi(\chi_{A_{n,1}}u) . \]
Thus
\begin{eqnarray*}
1 & = & \Phi(u) \\
  & = & \sumj \Phi(\chi_{A_{n,j}}u) \\
  & = & \ktwon \Phi(\chi_{A_{n,1}}u) .
\end{eqnarray*}
Therefore,
\begin{eqnarray*}
\phinj(y) & = &  \Phi(\chi_{A_{n,j}}y) \\
  & = & y_n(j)(\ktwon)^{1/p}\Phi(\chi_{A_{n,j}}u) \\
  & = & y_n(j)(\ktwon)^{1/p}\Phi(\chi_{A_{n,1}}u) \\
  & = & \frac{y_n(j)}{(\ktwon)^{1/q}} .
\end{eqnarray*}
Thus $P_ky = y$ for all $y \in Y_k$, as required.
\end{pf}

Propositions \ref{embed} and \ref{eachk} combine to give

\begin{prop}
For all $k \in{\rm\N}$, $L^{p,\infty}[0,k]$ is $q-$isomorphic to a
subspace of $X_k$ which is complemented in $X_k$ by a projection of
norm $\leq  q^2$.
\end{prop}

The following Theorem, the main goal of this section, now follows
easily.

\begin{thm} \label{first}
The space $(\sum^{\infty}_{k=1}\oplus
L^{p,\infty}[0,k])_{\ell^{\infty}}$ is isomorphic
to a complemented subspace of $(\sum^{\infty}_{k=1}\oplus
X_k)_{\ell^{\infty}}$ .
\end{thm}

\section{Embedding $X_1$ into $\seq$ complementably}

Choose a strictly increasing sequence of integers
$(m_n)^{\infty}_{n=0}$ so that $m_0 = 1$ and
\[ m_n \geq \sum^{n-1}_{j=1}2^jm_j \hspace{1em} \mbox{for all\ } n \geq
1. \]
Then choose a pairewise disjoint sequence of subsets of\ \,\N,
$(B_{n,j})^{2^n\ \infty}_{j=1 n=0}$, so that $|B_{n,j}| = m_n$, where
$|B|$ is the cardinality of the set $B$. Finally, let $w_{n,j} =
\chi_{B_{n,j}}$ for $1\leq j\leq 2^n, n\geq 0$. If $W$ is the
weak$^*$ closed subspace of $\seq$ generated by $(w_{n,j})$, i.e.,
\[ W = \{(a_i)\in \seq : (a_i) \mbox{\ is constant on each\ } B_{n,j}
\} , \]
then we will show that $X_1$ is isomorphic to $W$. Since $W$ is
complemented in $\seq$ by the conditional expectation operator with
respect to the $\sigma$-algebra generated by $(B_{n,j})$, we will have
proved

\begin{thm} \label{second}
$X_1$ embeds complementably into $\seq$.
\end{thm}

To show the isomorphism between $X_1$ and $W$, define
\[ Rx =
\bigvee^{\infty}_{n=0}\bigvee^{2^n}_{j=1}\frac{x_n(j)}{m^{1/p}_n}
w_{n,j}  \]
for all $x\in X_1$, where the suprema refer to the pointwise order on
the vector lattice of all real sequences.  We first show that $R$ maps
into $\seq$.  Fix $x \in X_1$.  Let $A$ be the set of all ordered
pairs $(n,j)$, such that $1\leq j\leq 2^n$, and
\begin{equation} \label{defA}
 |m^{-1/p}_nx_n(j)| < |m^{-1/p}_lx_l(i)|
\end{equation}
for some $l > n,\ 1 \leq i \leq 2^l$.
Since $\sup_j|m^{-1/p}_nx_n(j)| \leq m^{-1/p}_n\|x\| \to 0$ as $n\to
\infty$, for every $(n,j)\in A$, there exists $(l,i)\in A^c$
satisfying equation (\ref{defA}), with $l > n$.
Let $\Psi:A\to A^c$ be a choice function such that $\Psi(n,j) = (l,i)$
satisfies (\ref{defA}) with respect to $(n,j)$, and $l > n$. Let
\[ C_{l,i} = \bigcup_{(n,j)\in\Psi^{-1}\{(l,i)\}}B_{n,j} \]
for all $(l,i)\in A^c$.  Then
\begin{eqnarray*}
\sum_{(n,j)\in\Psi^{-1}\{(l,i)\}}\lv\frac{x_n(j)}{m^{1/p}_n}\rv w_{n,j}
 & < &
\lv\frac{x_l(i)}{m^{1/p}_l}\rv\sum_{(n,j)\in\Psi^{-1}\{(l,i)\}}w_{n,j}
\\
 & = & \lv\frac{x_l(i)}{m^{1/p}_l}\rv\chi_{C_{l,i}} .
\end{eqnarray*}
Note that
\begin{eqnarray*}
|C_{l,i}| & \leq & \sumn\sumjn|B_{n,j}| \\
 & = & \sumn\sumjn m_n \\
 & \leq & m_l
\end{eqnarray*}
by the choice of $(m_n)$.  It follows that the sequence
$\vee_{(n,j)\in A}m^{-1/p}_n|x_n(j)|w_{n,j}$
can be rearranged so that it is $\leq \vee_{(n,j)\notin
A}m^{-1/p}_n|x_n(j)|w_{n,j}$  in the
pointwise order.  Since the norm of $\seq$ is rearrangement invariant,
we have
\begin{equation} \label{parts}
  \lnorm\bigvee_{(n,j)\in A}m^{-1/p}_nx_n(j)w_{n,j}\rnorm \leq
  \lnorm\bigvee_{(n,j)\notin A}m^{-1/p}_nx_n(j)w_{n,j}\rnorm .
\end{equation}
On the other
hand, for all $n \geq 0$, let  $(x^*_n(j))^{2^n}_{j=1}$ denote the
 decreasing
rearrangement of $(|x_n(j)|)^{2^n}_{j=1}$, then
\begin{eqnarray*}
\lnorm\bigvee_{(n,j)\notin A}m^{-1/p}_nx_n(j)w_{n,j}\rnorm & \leq &
    q\lv\lv\lv\bigvee_{(n,j)\notin A}m^{-1/p}_nx_n(j)w_{n,j}\lv\lv\lv \\
 & \leq & \sup_{n,j}\frac{x^*_n(j)\lp
     j\!\cdot\!m_n+\sum^{n-1}_{i=1}2^im_i\rp^{1/p}}{m^{1/p}_n} \\
 & \leq & \sup_{n,j}\frac{x^*_n(j)\, ((j+1)m_n)^{1/p}}{m^{1/p}_n} \\
 & \leq & 2^{1/p}\sup_{n,j}j^{1/p}x^*_n(j) \\
 & = & 2^{1/p}\sup_n|||x_n||| \\
 & \leq & 2^{1/p}\sup_n\|x_n\| \\
 & = & 2^{1/p}\|x\| .
\end{eqnarray*}
Together with equation (\ref{parts}), this shows that $R$ is bounded
as a map into $\seq$.  On the other hand, for all $x\in X_1$,
\begin{eqnarray*}
 \|Rx\| & \geq &
  \sup_n\lnorm\sum^{2^n}_{j=1}\frac{x_n(j)}{m^{1/p}_n}w_{n,j}\rnorm \\
 & = & \|x\|.
\end{eqnarray*}
Therefore, $R$ is an embedding, as claimed.

\section{Proof of the main theorem}

\begin{prop} \label{last}
$\rline$ embeds complementably into
$\lp\sum^{\infty}_{k=1}\oplus L^{p,\infty}[0,k]\rp_{\ell^{\infty}}$.
\end{prop}

\begin{pf}
This is rather straight forward.  Define $S:\rline\to
\lp\sum^{\infty}_{k=1}\oplus L^{p,\infty}[0,k]\rp_{\ell^{\infty}}$ by
\[ Sf = (f\chi_{[0,1]},f\chi_{[0,2]}, \ldots ) \]
for all $f\in \rline$.  Clearly, $S$ is an isometric embedding. Now
choose a free ultrafilter \,${\cal U}$ on \N, and regard
$L^{p,\infty}[0,k]$ as the subspace of $\rline$ consisting of all
functions supported on $[0,k]$.  Since $\rline$ is the dual of
$L^{q,1}[0,\infty)$, its unit ball is weak$^*$ compact.  Given $g =
(g_k) \in \lp\sum^{\infty}_{k=1}\oplus
L^{p,\infty}[0,k]\rp_{\ell^{\infty}}$, let
\[ Qg = (w^*)\lim_{k\to{\cal U}}g_k\ . \]
It is easy to see that $Q$ is bounded as a map into $\rline$, and that
$Q\circ S$ is the identity on $\rline$.  This proves the Proposition.
\end{pf}

The proof of Theorem \ref{main} now follows as in the discussion in
\S\ref{theorem}, using Theorems \ref{first}, \ref{second}, and
 Proposition
\ref{last} above.

\flushleft
\vspace{.5in}
Department of Mathematics\\National University of Singapore\\
Singapore 0511\\ e-mail(bitnet) : matlhh@nusvm

\end{document}